%%%%%%%%%%%%%%%%%%%%%%%%%%%%%%%%%%%%%%%%%%%%%%%%%%%%%%%%%%%%%%%%%%%%%%%%%%%%%%

\documentstyle{amsppt}
\nologo
%\NoPageNumber
%
% ------ Macros ------
%
\font\b=cmr10 scaled \magstep4
\def\bigzerou{\smash{\kern-25pt\lower1.7ex\hbox{\b 0}}}
\hsize=360pt
\vsize=500pt
\hbadness=5000
\tolerance=1000
\NoRunningHeads

\def\z{\zeta}
\def\Q{\Bbb Q}

\def\al{\alpha}

\def\Fp{\Bbb F_p}
\def\F5{\Bbb F_5}

\def\Z{\Bbb Z}
\def\Zp{\Bbb Z_{(p)}}

\def\CC{\Cal C}

\def\no{\noindent}

\def\vp{\varphi}
\def\i{\infty}
\def\ti{\times}

\def\Si{\Sigma}

\def\Ga{\Gamma}

\def\ot{\otimes}
\def\part{\partial}

\magnification=\magstep1
%\pageno=27
\topmatter
\title On $K_1$ of a self-product of a curve\endtitle
\author
Kenichiro Kimura
\endauthor
\footnote{The author was partially supported by JSPS Postdoctral Fellowships
for Research Abroad for 2000-2001.}

\affil Institute of Mathematics, University of Tsukuba,
Tsukuba,305-8571, Japan 
email: kimurak\@math.tsukuba.ac.jp\endaffil

\abstract{We present elements of $H^1(C\ti C,\,\Cal K_2)$ for certain
specific curves $C$. 
The image of the element under the boundary map 
arising from the localization sequence of $K$-theory
is the graph of frobenius
endomorphism of the reduction of the curve modulo 3.}
\endabstract
\endtopmatter
\vskip -3ex
\noindent\hskip 2em

\vskip 3ex
\CenteredTagsOnSplits
% ------  Document  ------
%
\document

\subhead{\S 1 Introduction}\endsubhead
 Let $X$ be a projective smooth variety over
$\Q$. Let S be a finite set of primes such that $X$ has 
a projective smooth model $\Cal X$ over $U:=\text{spec}\Z[\frac1S]$.
Then from the works of Bloch(\cite{Bl1}) and Sherman, 
we have the following exact localization sequence in algebraic 
$K$-theory, using Gersten's conjecture for $\Cal K_2$
in a mixed characteristic setting that was proven by Bloch(\cite{Bl2})

$$H^1(\Cal X,\,\Cal K_2)\to H^1(X,\,\Cal K_2)
\overset{\part}\to{\to} \underset{p\in U}\to\oplus Pic(X_p)
\to CH^2(\Cal X)\to CH^2(X)\to 0$$ where $X_p$ is the fiber of
$\Cal X$ at the prime $p$.

A conjecture of Beilinson on the special values of $L$-function
and Tate conjecture tell us that the cokernel of $\part$
is torsion. See the remark after Theorem 2.5 in \cite{La}
for an account of this.
This means there should be enough elements
in $H^1(X,\,\Cal K_2)$ so that any element of $Pic(X_p)\ot\Q$ is in the
image of the map $\part\ot\Q$. The purpose of this paper
is to give a potential method to construct elements in
$H^1(X,\,\Cal K_2)$ in the case where
$X$ is the product of a curve.

\no In this case a typical element of $Pic(X_p)$ which usually
does not lift to the generic fiber is the graph
of frobenius endomorphism of the fiber of the 
curve.  For a certain curve $C$ 
of special type we construct 
a divisor $D\subset C\ti C$ such that the reduction
$D_p$ of $D$ modulo $p$ 
contains the graph of frobenius endomorphism $\Ga_{\vp}$
as an irreducible component(Lemma 2.1). 

Next task is to find a nice function $g$ on $D$
such that $g$ vanishes on $\Ga_{\vp}$ and not on 
other components of $D_p$, and that one can 
modify $(g,D)$ to an element $\Si\in H^1(X,\,\Cal K_2)$
by adding some elements. In this we have succeeded
for the curve $C$ and for the prime 3. This is the
main result of this paper (Theorem 2.3).

Although in this paper we are able to construct 
elements only for curves of very specific type,
 the idea of construction can be
applied to more general curves, and we hope it gives
a clue to construct more elements in more general 
situation.

When $X$ is the product of a modular curve there is a 
beautiful construction of elements in $H^1(X,\,\Cal K_2)$
by Flach(\cite{Fl}) and Mildenhall(\cite{Mi}).
 However their
elements are still not enough to show that coker($\part$)
is torsion for the products of general modular curves.

The author thanks S. Bloch and R. Sreekantan for discussions.
 This work
was carried out while the author stayed at Department of Mathematics,
University of Chicago. He is thankful to the people at Chicago 
for their hospitality.

\subhead{\S 2 Construction of an element in
$H^1(X,\,\Cal K_2)$}\endsubhead

 Let $C$ be a curve over $\Q$ defined as the normalization
of the projective closure of the affine curve 
$$y^m=t^n+1.$$ Here $m$ and $n$ are integers such that
$n\geq 5,\,\,(n,6)=1,\,\, (m,3)=1$ and $(m,n)=1$. 
We denote its point at infinity by $\i$.
For a prime $p\nmid mn$, the curve $C$ naturally
extends to a projective smooth model $\Cal C$ over 
spec $\Zp$. We denote the special fiber of $\Cal C$ by
$C_p$. Let  
$$\pi: \,\,\CC\to \Bbb P^1_{\Zp}$$ be a map given by
$(y,\,t)\mapsto (t)$. Let $F:\,\Bbb P^1_{\Zp}\to \Bbb P^1_{\Zp}$
be a morphism given by $t\mapsto t^p.$

\no Let $f:=F\circ \pi$ and $\Ga_f\subset 
\CC\ti_{\Zp} \Bbb P^1$ be its graph. Let $\Cal D=(id\ti \pi)^*
\Ga_f \subset \CC\ti\CC$. We denote the generic fiber of $\Cal D$
by $D$ and the special fiber of $\Cal D$ by $D_p$.

\proclaim{Lemma 2.1} The intersection 
$$\Cal D\cap(\CC\ti \i \cup \i\ti \CC)$$ is supported 
on $\i\ti \i.$\endproclaim

\demo{proof} If  $x\ti \i\in D_p\,\,(\text{resp.}\quad
x\ti \i\in D)$ for a closed point $x$ of $C_p\,\,(\text{resp. of}\,\,
C)$
then $x\ti \pi(\i)
=x\ti \i \in \Ga_f$ so $ F\circ \pi(x)=\i$, hence 
$\pi(x)=\i,\,\,x=\i.$

\no Similarly, if
$\i\ti x\in D_p\,\,(\text{resp.}
\in D) \Rightarrow \i\ti\pi(x)\in \Ga_f,
\quad \pi(x)=f(\i)=\i, x=\i.$ \qed \enddemo

\proclaim{Lemma 2.2} Let $\kappa$ be the field generated 
by a primitive $m$-th root of unity $\xi$ over $\Fp$. Then the divisor

$$D_p\ot_{\Fp} \kappa=\sum_{j=0}^{m-1}\Ga_{\xi^j\circ \vp}.$$
Here $\vp$ is the frobenius endomorphism 
of $C_p$ and $\xi$ is the automorphism
of $\CC$ given by $y\mapsto \xi y.$\endproclaim
\demo{proof}
 Let $\Cal V:=(\CC\backslash \{\i\})\ti_{\Zp} (\CC \backslash \{\i\})$. 
We denote the generic fiber of $\Cal V$ by $V$ and its specical
fiber by $V_p$.

\no  By lemma 2.1 it suffices to look at $D_p\cap V_p$. 
On $V_p$, $D_p$ is defined
by the equation $t_2-t_1^p=0$. Here $(y_1,t_1)\ti(y_2,t_2)$ is a set of
coordinates on $V_p$. One can see that $y_2^m-(y_1^p)^m\in 
(t_2-t_1^p)$. \qed \enddemo

\remark{Remark}  This construction of the divisor $D$
can be applied in general to curves which are cyclic
cover of $\Bbb P^1$.
\endremark

In the following we specify $p$ to 3. When $m=2$ and $n=5$ 
the curve
$C_3$ over $\Bbb F_3$ is isomorphic to the curve
$t^5=y(y-1)$. In this case the characteristic 
polynomial of the action of $\vp$ 
on the Tate module of Jacobian of 
$C_3$ is $T^4+9$(\cite{Gr-Ro} Lemma1.1). Since $End(Jac(C))=\Z$,
$\Ga_\vp$ does not lift to the generic fiber. 

\proclaim{Theorem 2.3} There is an element $\Si \in H^1(C\ti C,\,\Cal K_2)
\ot \Q$
such that the factor of 
$$\part(\Si)\,\, at\,\,p\,=\cases & \Ga_\vp\quad p=3\\
                    & 0  \quad p\in \text{spec}\Z[\frac{1}{mn}],\,\,p\neq 3.
  \endcases $$\endproclaim
\demo{proof}  Consider the function $g:=y_2-y_1^3$ on $D$.
We calculate div$(g)\cap V$.

\no $y_2 -y_1^3=0\Rightarrow
t_2^n+1-(t_1^n+1)^3=0.$ 

\no We write $T=t_1^n$. Then it follows
$T^3+1-(T+1)^3=0,\,-3(T^2+T)=0,\,T=0\,\,\text{or}\,\,-1$.

\no So $t_1=0,-\z^\al,\,\al=0,\cdots , n-1.$ 

\no Here $\z$ is a primitive $n$-th root of unity.

\no Hence div$(g)$
is suppoted on the set $\{( \xi^j,0)\ti(\xi^{3j},0)|j=0,\cdots,m-1\}
\cup\{(0,-\z^\al)\ti(0,-\z^{3\al})|\al=0,\cdots , n-1\}\cup \{\i\}.$

\no Here $\xi$ is a primitive $m$-th root of unity.

\no Let $\Cal S$ be the set of points on $C$ given by

$$\Cal S=\{(y,\,t)\in C\backslash \{\i\}
| y=0\,\,\text{or}\,\,t=0\}\cup\{\i\}.$$
We see that the support of div$(g)\subset \Cal S\ti\Cal S.$

\no For any point $x\in \Cal S$ the divisor $(x)-(\i)\in Pic^0C$
is torsion, so we can construct an element $\Si \in H^1(C\ti C,\,\Cal K_2)
\ot \Q$
of the form 
$$(g,D)+\sum_{x_i\in \Cal S}(f_i,\, C\ti x_i)+(f'_i,\, x_i\ti C).$$
Since $g=0$ at the generic point of $\Ga_\vp$ and $g\neq 0$
at the generic points of other irreducible components
of $D_3$, we see that the factor of
$\part(\Si)$ at $p=3$ is equal to $\Ga_\vp.$

\proclaim{Lemma 2.4} For any prime $l\in \text{spec}\Z[\frac{1}{mn}]$
 other than 3, the reduction
$D_l$ of the divisor $\Cal D$ mod $l$ is irreducible.\endproclaim
\demo{proof} One can show that
$$D_l\cap (\i\ti C_l \cup C_l\ti \i)$$ is supported on
$\i\ti \i$ in the same way as the proof of Lemma 2.1.
So it suffices to show that
$D_l\cap (C_l\backslash \{\i\})\ti(C_l\backslash \{\i\})$ is irreducible.
 On $ (C_l\backslash \{\i\})\ti(C_l\backslash \{\i\})$, $D_l$ is the scheme

$$\text{spec}\,\,\Bbb F_l[y_1,t_1,y_2]/(y_1^m-t_1^n-1,\,\,
y_2^m-(y_1^m-1)^3-1).$$

\no So $D_l$ is reducible if and only if 
$y_2^m-(y_1^m-1)^3-1$ is reducible. There is a point 
$s\in C_l(\overline{\Bbb F_l})$ at which the function $(y_1^m-1)^3+1$
vanishes to order 1. It follows that $y_2^m-(y_1^m-1)^3-1$
is irreducible from the following well known fact:

{\it Let $\Cal O$ be a discrete valuation ring and denote its
maximal ideal by $\frak m$. Let
$f(x)=x^k+a_{k-1}x^{k-1}+\cdots +a_0\in \Cal O[x]$ be 
a polynomial. If $a_j\in \frak m$ for $j=0,\cdots,k-1$
and $a_0$ generates $\frak m$, then $f(x)$ is irreducible.} \qed

\enddemo

\no The above argument and Lemma 2.4 show the existance of $\Si$ 
as stated in the theorem. \qed \enddemo


\Refs

\widestnumber\key{Gr-Ro}

\ref
\key Bl1
\book Lectures on Algebraic Cycles
\by Bloch, S.
\publ Duke University Press
\yr 1980
\endref


\ref
\key Bl2
\paper ``A note on Gersten's conjecture in the mixed
characteristic case''in Applications of Algebraic
$K$-theory to Algebraic Geometry and Number Theory,
Part I, II
\by Bloch, S.
\jour Contemp. Math.
\vol 55
%\pages 75-88
\yr 1986
\publ Amer. Math. Soc. Providence
\endref

\ref
i
\key Fl
\paper A finiteness theorem for the symmetric square
of an elliptic curve
\by  Flach, M.
\jour Invent. Math
\yr 1992
\pages 307-327
\vol 109
\endref



\ref
i
\key Gr-Ro
\paper Some results on the Mordell-Weil group of 
the Jacobian of the Fermat curve
\by  Gross, B., Rohrlich, D. 
\jour Invent. Math
\yr 1978
\pages 201-224
\vol 44
\endref




\ref
\key La
\paper Zero cycles on Hibert-Blumenthal surfaces
\jour Duke Math. J
\vol 103
\yr 2000
\pages 131-163
\by Langer, A.
\endref




\ref
\key Mi
\paper Cycles in a product of elliptic curves, and a group
analogous to the class group
\jour Duke Math. J
\vol 67
\yr 1992
\pages 387-406
\by Mildenhall, S.
\endref

\endRefs
\enddocument